\documentclass{article}                   

\usepackage{graphicx}
\usepackage{mathptmx}     
\usepackage{cite}

\usepackage{latexsym}
\usepackage{amssymb}
\usepackage{amsmath}

\newcommand{\bN}{\mathbb{N}}

\newcommand{\cE}{\mathcal{E}}

\newcommand{\cG}{\mathcal{G}}
\newtheorem{theorem}{Theorem}[section]
\newtheorem{lemma}{Lemma}
\begin{document}

\sloppy
 
\title{\itshape{Exact Distribution of the Sample Variance from a Gamma Parent   
Distribution}}
  
\author{ T. Royen \\
Fachhochschule Bingen, University of Applied Sciences,\\
Berlinstrasse 109, D--55411 Bingen, Germany\\
e-mail: royen\symbol{64}fh-bingen.de} 
\date{}
\maketitle

\begin{abstract}
Several representations of the exact cdf of the sum of squares of $n$ independent identically gamma--distributed random variables $X_i$ are given, in particular by a series of gamma distribution functions. Using a characterization of the gamma distribution by Laha, an expansion of the exact distribution of the sample variance is derived by a Taylor series approach with the former distribution as its leading term. In particular for integer orders $\alpha$ some further series are provided, including a convex combination of gamma distributions for $\alpha = 1$ and nearly of this type for $\alpha > 1$. Furthermore, some representations of the distribution of the angle $\Phi$  between $(X_1,..., X_n)$ and $(1,...,1)$ are given by orthogonal series. All these series are based on the same sequence of easily computed moments of $\cos (\Phi)$.\\


\end{abstract}

\section{Introduction}
\label{intro}
The distribution of the sample variance $s^2$ in non--normal cases has attracted sporadic attention during the last eight decades. Early investigations, concerning a gamma parent distribution can be traced at least to Craig (1929)  and Pearson (1929), who used moment approximations, later investigated more thoroughly by Bowman and Shenton(1983). For general distributions various approximations were proposed by Box (1953), Roy and Tiku (1962), Tan and Wong (1977) and by Mudholkar and Trivedi (1981). The latter authors recommend transformations of the Wilson--Hilferty type. Exact results seem to be available only for a mixture of two normal distributions (Hyrenius 1950, Mudholkar and Trivedi 1981). See also chapter~4.7 in Mathai and Provost (1992). The interest in the distribution of $s^2$ or $1/s$ was also stimulated by investigations on the estimation of process capability indices, see e.g. Pearn, Kotz and Johnson (1992) and Pearn and Kotz (2006).

Apart of its intrinsic theoretical value, an analytical representation of the cdf of $s^2$ enables more rapid and more accurate evaluations than Monte Carlo methods. Besides, it provides a base to investigate the accuracy of the different proposed approximations.

Now let $X_1,\ldots,X_n$ be a random sample with mean $\bar{X}$ and sample variance\linebreak \mbox{$S^2 = (n-1)^{-1} \sum^{n}_{i=1} (X_i - \bar{X})^2$} from a gamma distribution with pdf
\begin{equation}
\label{eq:1a}   
g_\alpha(x) = x^{\alpha -1} \exp(-x)/\Gamma (\alpha),
\end{equation}
and cdf
\[
G_\alpha (x) = \int^{x}_{0} g_\alpha (\xi)d\xi .
\]
With
\begin{equation}
\label{eq:1b}    
Y:= \sum^{n}_{i=1}X_i, \quad Z:= \sum^{n}_{i=1}X^2_i, \quad U^2:= ZY^{-2},
\end{equation}
we have
\begin{equation}
\label{eq:1c}   
(n-1) S^2 = (U^2 - 1/n) Y^2, \quad  1/\sqrt n \leq U \leq 1.
\end{equation}
By Laha (1954) --- see also Lukacs (1955) --- the mutual independence of $\overline{X}$ and the sample coefficient of variation  $S/\overline{X} = \sqrt n((U^2 - 1/n)/(1-1/n))^{1/2}$ was shown to characterize the family of gamma distributions indexed by $\alpha$. This is the key to derive the distribution of $S^2$ essentially from the distribution of $Z$. At first, it follows with the angle $\Phi$ between $(X_1,\ldots,X_n)$ and $(1,\ldots,1)$ that 
\begin{eqnarray}
\label{eq:1d}     
H_{\alpha,n}(r) &:=& \Pr \{Z\leq r^2\} \ = \ {\cal E}_U (G_{\alpha n}(r/U))\nonumber \\
&&\\
&=& {\cal E}_\Phi (G_{\alpha n} (r \sqrt n \cos \Phi)),\nonumber
\end{eqnarray}
\begin{eqnarray}
\label{eq:1e}     
\Pr \{ \sqrt{n-1} S \leq r\} &=& {\cal E}_U (G_{\alpha n} (r(U^2-1/n)^{-1/2})\nonumber\\
&&\\
&=& {\cal E}_\Phi (G_{\alpha n} (r\sqrt{n} \cot \Phi)).\nonumber
\end{eqnarray}
With the cdf $F_{\alpha,n}$ of $U$ it follows also
\begin{eqnarray}
\label{eq:1f}   
H_{\alpha,n}(r) &=& {\cal E}_Y (F_{\alpha,n} (r/Y)) \nonumber\\
\\
&=& r^{\alpha n}(\Gamma(\alpha n))^{-1} \int^{\sqrt n}_{0} x^{\alpha n-1} F_{\alpha,n} (x^{-1}) e^{-rx} dx.\nonumber
\end{eqnarray}
Thus
\begin{equation}
\label{eq:1g}   
h^*_{\alpha,n} (r) := r^{-\alpha n} H_{\alpha,n}(r)
\end{equation}
is the Laplace transform (L.t.) of
\begin{equation}
\label{eq:1h}    
f^*_{\alpha,n}(x) := (\Gamma(\alpha n))^{-1} x^{\alpha n-1} F_{\alpha,n}(x^{-1}), \quad 0 < x \leq \sqrt n.
\end{equation}
A power series for $h_{\alpha,n}^*$ -- and consequently the cdf $H_{\alpha,n}$ -- is derived in  (\ref{eq:2d}) $\ldots$ (\ref{eq:2h}) of the follwing section. Due to  (\ref{eq:1g}),  (\ref{eq:1h}) the coefficients of this series determine the moments of $f_{\alpha,n}^*$  and therefore $F_{\alpha,n}$ and the cdf in  (\ref{eq:1e}).

For a cdf a probability mixture representation is more appealing than a power series. Such a mixture of gamma distributions for $H_{\alpha,n}$ is found in (\ref{eq:2m}) in section \ref{sec:2}.

Then the reader might go directly to theorem \ref{theorem:4.1}, which is the main result concerning the actual computation of the cdf of $S^2$. By means of a few lines of code of a computer algebra system some tables of this cdf with at least eight correct digits have been computed for some values of $\alpha$ and  $n$.

Representations for $f_{\alpha,n}^*$ by an orthogonal series with Legendre polynomials and by a Fourier sinus series are derived too in section \ref{sec:2}. Because of the relations (\ref{eq:1g}),  (\ref{eq:1h}) these two series provide two further representations of $h_{\alpha,n}^*$ or $H_{\alpha,n}$. The corresponding series for $F_{\alpha,n}(u)$ in (\ref{eq:3a}) and the orthogonal series derived from  (\ref{eq:3b}) can be used within (\ref{eq:1e}).  However, the application of theorem \ref{theorem:4.1} is simpler and more accurate.

In spite of the numerical use of the double series in  theorem \ref{theorem:4.1} it would be theoretically more satisfying to have a single alternating power series or even a probability mixture representation for the cdf of $S^2$. The latter is accomplished exactly for the exponential case $(\alpha = 1)$ and nearly for integer $\alpha >1$ by theorem \ref{theorem:4.2}, which contains also an alternating power series. However, for numerical purposes, theorem \ref{theorem:4.2} is (at present) not considered as a competitor to theorem \ref{theorem:4.1} since the computation of the required coefficients is more cumbersome. The method is explained in section \ref{sec:3}. It is based on a representation of the cdf of $\tan \Phi$ by a polynomial on a certain section $[0,\varphi_n]$ of the domain of $\Phi$. This formula is found in theorem \ref{theorem:3.1}. The proofs of all the theorems are given in the appendix. 

Throughout the paper formulas from the handbook of mathematical functions by Abramowitz and Stegun  are cited by A.S. and their number. The symbol $\sum_{(k)}$ stands for summation over all possible decompositions $k = k_1+\ldots +k_n$ with non--negative integers $k_j$. Moments of positive random variables occur also with non-integer exponents and for defective distributions.

\section{The distribution of the sum of squares}
\label{sec:2}
The L.t. of $Z_1 = X^2_1$ is given by
\begin{eqnarray}
\label{eq:2a}  
\psi_\alpha(t) &=& (2\Gamma (\alpha))^{-1} \int_0^\infty z^{\alpha/2-1}\exp(-\sqrt z - tz)dz\nonumber \\
&=& (2\Gamma (\alpha))^{-1} \sum^{\infty}_{k=0} (-1)^k \ \frac{\Gamma((\alpha+k)/2)}{k!} \ t^{-(\alpha + k)/2} \\
&=& (2\Gamma (\alpha))^{-1}t^{-\alpha/2} \exp\left(\frac 1{4t}\right)\nonumber \\
&&\cdot \left( \Gamma \left(\frac \alpha 2\right) M\Bigl( \frac{1-\alpha}2, \frac 12; -\frac 1{4t}\Bigr)
 - \frac 1{\sqrt t} \ \Gamma \Bigl(\frac{\alpha+1}2\Bigr) M\Bigl( \frac{2-\alpha}2, \frac 32; - \frac 1{4t}\Bigr) \right)\nonumber
\end{eqnarray}
with Kummer's confluent hypergeometric function $M$. In particular
\begin{equation}
\label{eq:2b}   
\psi_1(t) = \sqrt{\frac \pi t} \ \exp \left( \frac 1 {4t}\right) \Phi \left( - \frac 1{\sqrt{2t}}\right)
\end{equation}
with the cdf $\Phi(x) = (1 + {\mathrm{erf}} (x/\sqrt 2))/2$ of the standard normal distribution. For integer values $(1-\alpha)/2$ or $(2-\alpha)/2$ Kummer's $M$ is given by Hermite polynomials \linebreak  (A.S. 13.6.17/18).

The cdf $H_{\alpha,n}(r) = \Pr \{Z \leq r^2\}$ can be obtained by the Fourier inversion formula, but some further representations are useful. If $\psi_\alpha(t)$ is written as $t^{-\alpha/2} \beta_\alpha (t^{-1/2})$,  then
\begin{equation}
\label{eq:2c}   
(\psi_\alpha(t))^n = t^{-\alpha n/2} \sum^{\infty}_{k=0} \beta_{\alpha,n,k} \ t^{-k/2}
\end{equation}
with 
\begin{eqnarray}
\label{eq:2d}   
\beta_{\alpha,n,k} &=& (-1)^k (2\Gamma(\alpha))^{-n} \sum_{(k)} \prod^{n}_{i=1} \ \frac{\Gamma ((\alpha + k_i)/2)}{k_i!}\nonumber\\
&=& \frac 1\pi \int^{\pi}_{0} {\mathrm{Re}} \left\{ (\beta_\alpha (y^{-1}))^n y^k\right\} d\varphi, \quad y = \rho e^{i\varphi}. \\[1ex]
&& \rho > 0, \quad - \pi < \varphi \leq \pi, \quad k \in \bN_0. \nonumber
\end{eqnarray}
The parameter $\rho$ was inserted here only for numerical considerations. The $\beta_{\alpha,n,k}$ are more quickly computed recursively for $n = 1,2,4,8,\ldots$ than by the integrals.

Laplace inversion in (\ref{eq:2c}) implies 
\begin{equation}
\label{eq:2e}   
H_{\alpha,n}(r) = r^{\alpha n} \sum^{\infty}_{k=0}\ \frac{\beta_{\alpha,n,k}}{\Gamma(1+(\alpha n + k)/2)} \ r^k.
\end{equation}
By the series expansion of $G_{\alpha n}$ in (\ref{eq:1d}) we obtain with the moments
\begin{equation}
\label{eq:2f}   
\mu_{\alpha,n,k} := \cE \left((\sqrt n \ \cos \Phi)^{\alpha n +k}\right) = n^{(\alpha n +k)/2} \gamma_{\alpha,n,k}
\end{equation}
\begin{equation}
\label{eq:2g}   
H_{\alpha,n}(r) = \frac {r^{\alpha n}}{\Gamma(\alpha n)}\ \sum^{\infty}_{k=0} \ \frac{\mu_{\alpha,n,k}}{\alpha n + k} \ \frac {(-r)^k}{k!}
\end{equation}
and by comparison with (\ref{eq:2e}) the relation 
\begin{equation}
\label{eq:2h}   
\frac{(-1)^k}{k!} \mu_{\alpha,n,k} = 2\Gamma(\alpha n)\frac{\beta_{\alpha,n,k}}{\Gamma((\alpha n + k)/2)} .
\end{equation}
Also the following representations of the distribution functions of $\sum^n_{i=1} X^2_i, S^2, U$ and $\tan \Phi$ are based essentially on the sequence $(\mu_{\alpha,n,k})$ or equivalently $(\gamma_{\alpha,n,,k})$.

A probability mixture representation for $H_{\alpha,n}$ by gamma distribution functions is obtained as follows: 
With any scale factor $\lambda$ let be 
\begin{equation}
\label{eq:2k}  
\mu_{\alpha,n,k,\lambda} = \lambda^{-(\alpha n + k)} \mu_{\alpha,n,k}\ .
\end{equation}
Then, multiplying
\[
\frac d{dr} H_{\alpha,n} (r) = \lambda \frac{(\lambda r)^{\alpha n-1}}{\Gamma(\alpha n)} \ \sum^{\infty}_{k=0} \mu_{\alpha,n,k, \lambda} (-\lambda r)^k/k!
\]
by $e^{-\lambda r} \sum^{\infty}_{\ell=0} (\lambda r)^\ell/\ell!$ we obtain
\begin{equation}
\label{eq:2ka}   
\lambda g_{\alpha n} (\lambda r) \sum^{\infty}_{k=0} \delta_{\alpha,n,k,\lambda} (\lambda r)^k/k!
\end{equation}
with the alternating $k^{\mathrm{th}}$ order differences 
\begin{equation}
\label{eq:2l}    
\delta_{\alpha,n,k,\lambda} = (-\Delta)^k \mu_{\alpha,n,0,\lambda} = \sum^{k}_{j=0} (-1)^j \binom{k}{j} \mu_{\alpha,n,j,\lambda}
\end{equation}
and consequently 
\begin{equation}
\label{eq:2m}   
H_{\alpha,n} (r) = \sum^{\infty}_{k=0}\delta_{\alpha,n,k,\lambda} \binom{\alpha n + k-1}{k} G_{\alpha n + k} (\lambda r).
\end{equation}
In particular with $\lambda = \sqrt n$ this is a probability mixture since \newline 
$\delta_{\alpha,n,k,\sqrt n} = \cE((\cos \Phi)^{\alpha n} (1-\cos \Phi)^k) > 0$ and
\begin{equation}
\label{eq:2ma}   
\sum^{\infty}_{k=0}  \binom{\alpha n + k-1}{k}\delta_{\alpha,n,k,\sqrt n} =1
\end{equation}
by virtue of the binomial series with $0 \leq 1 - \cos \Phi \leq 1 - 1/\sqrt n$. However, for numerical approximations different values of $\lambda$ should be more suitable, e.g.
\begin{equation}
\label{eq:2n}   
\lambda = \left( \cE \left( U^{-\alpha n}\right) \right)^{1/(\alpha n)} = \mu^{1/(\alpha n)}_{\alpha,n,0} = \left( \frac{2\Gamma(\alpha n)}{(2\Gamma(\alpha))^n}\frac{(\Gamma(\alpha/2))^n}{\Gamma(\alpha n/2)}\right)^{1/(\alpha n)}
\end{equation}
with $\mu_{\alpha,n,0,\lambda} =1$.\\

Besides, two further representations of $H_{\alpha,n}$ are obtained from the following orthogonal series for $f^*_{\alpha n}$ from (\ref{eq:1g}), (\ref{eq:1h}). With 
\begin{equation}
\label{eq:2o}     
\int^{\sqrt n}_{0} x^k  f^*_{\alpha n}(x) dx = \frac{1}{\Gamma(\alpha n)} \frac{\mu_{\alpha,n,k}}{\alpha n+k} \ ,
\end{equation}
implied by (\ref{eq:2f}), (\ref{eq:2g}), and the shifted Legendre polynomials
\begin{equation}
\label{eq:2p}  
P^*_k(y) = \sum^{k}_{j=0} p^*_{k,j} y^j, \ 0 \leq y \leq 1, \ \int^{1}_{0} (P^*_k(y))^2 dy = \frac{1}{2k+1} \ ,
\end{equation}
(A.S. 22.2.11), it follows 
\begin{eqnarray}
\label{eq:2q}    
f^*_{\alpha,n} (x) &=& \frac 1{\sqrt n} \sum^{\infty}_{k=0} c_{\alpha,n,k} P^*_k (x/\sqrt n) \ ,\nonumber \\[-1ex]
&& \\
c_{\alpha,n,k} &=& (2k +1) \frac{n^{\alpha n/2}}{\Gamma(\alpha n)} \sum^{k}_{j=0} p^*_{kj} \frac{\gamma_{\alpha,n,j}}{\alpha n +j}\ . \nonumber
\end{eqnarray}
Thus, with (A.S.~11.4.26) and the modified spherical Bessel functions (A.S.~10.2.2) we find  
\begin{equation}
\label{eq:2r}   
H_{\alpha,n}(r) = r^{\alpha n} \exp \left( - \frac 12 r \sqrt n\right)\sum^{\infty}_{k=0} (-1)^k c_{\alpha,n,k} \left( \frac \pi{r \sqrt n}\right)^{1/2} I_{k+1/2} \left( \frac 12 r \sqrt n\right) . 
\end{equation}
Consequently, the Fourier transform of $f^*_{\alpha,n}$ is representable by
\begin{equation}
\label{eq:2s}   
h^*_{\alpha,n}(-it) =\exp \left(\frac i2 \sqrt n t\right)\sum^{\infty}_{k=0} i^k c_{\alpha,n,k}\, \boldmath{j}_k \left( \frac 12 \sqrt n t\right)
\end{equation}
with the spherical Bessel functions $j_k$ from (A.S.~10.1). For absolutely large real $t$ this series is numerically more suitable than the power series from (\ref{eq:2g}). Such values are required as coefficients in the following Fourier sinus expansion for $f_{\alpha,n}^*$. \\

\begin{equation}
\label{eq:2t}    
f^*_{\alpha,n}(x) = \sum^{\infty}_{m=1} b_{\alpha,n,m} \sin\left( \frac{m\pi x}{\sqrt n} \right), \ 0 \leq x \leq \sqrt n\ ,
\end{equation}
with 
\begin{eqnarray}
\label{eq:2u}    
b_{\alpha,n,m} &=& \frac 2{\sqrt n} {\mathrm{Im}} \left( h^*_{\alpha,n} \left( -i \frac{m \pi}{\sqrt n}\right)\right)\nonumber \\[-1ex]
&&\\
&=& \frac 2{\sqrt n}    \sum_{k\geq 0 \atop m+k \ {\mathrm{odd}}} (-1)^{(m+k-1)/2} c_{\alpha,n,k} \,\boldmath{j}_k (m\pi/2) \nonumber\ .
\end{eqnarray}
Because of lemma~\ref{lemma:2} the function $f^*_{\alpha n}$ has a square integrable derivative at least for $\min (\alpha n,n)>2$ which entails the absolute uniform convergence of the orthogonal series in (\ref{eq:2q}) and (\ref{eq:2t}). Inserting (\ref{eq:2t}) into (\ref{eq:1f}) leads to 
\begin{equation}
\label{eq:2v}   
H_{\alpha, n} (r) = r^{\alpha n} \sum^{\infty}_{m=1} \frac{m\pi \sqrt n}{nr^2 + m^2\pi^2} \ b_{\alpha,n,m} \left(1-(-1)^m e^{-r \sqrt n}\right) .
\end{equation}
Finally, with the functions 
\begin{eqnarray*}
\label{eq:2vi}   
\lefteqn{\cG_\beta (\zeta) := \sum ^{\infty}_{k=0} \ \frac{\zeta^k}{\Gamma(\beta+k/2)}\nonumber}\\
&=& (\Gamma(\beta))^{-1} M(1,\beta;\zeta^2) + (\Gamma(\beta + 1/2))^{-1} \zeta M(1,\beta+1/2;\zeta^2), \ \beta >0,\\
&=& \zeta^{2(1-\beta)} e^{\zeta^2}\left(G_{\beta-1}(\zeta^2) + G_{\beta - 1/2}(\zeta^2)\right), \quad {\mathrm{if}} \quad \beta > 1, \nonumber
\end{eqnarray*}
the following integral representation can be derived from (\ref{eq:2c}), (\ref{eq:2d}), (\ref{eq:2e}):
\begin{eqnarray*}
\label{eq:2vj}    
H_{\alpha,n}(r) &=& \frac{r^{\alpha n}}{\pi} \int^{\pi}_{0} {\mathrm{Re}} \left\{\left(\beta_\alpha(y^{-1})\right)^n \cG_{1+\alpha n/2} (ry)\right\} d\varphi,\\
y &=& \rho e^{i\varphi}, \quad \rho > 0.\nonumber
\end{eqnarray*}

\section{The distribution of the angle ${\boldsymbol{\Phi}}$}
\label{sec:3}
Two representations of the cdf $F_{\alpha,n}$ of $U=\left( \sqrt n \cos \Phi\right)^{-1}$ follow directly from (\ref{eq:2q}) and (\ref{eq:2t}).
\begin{eqnarray}
\label{eq:3a}   
F_{\alpha,n}(u) &=& \Gamma(\alpha n) u^{\alpha n - 1} \frac 1{\sqrt n} \sum^\infty_{k=0} c_{\alpha,n,k} P^*_k\left( 1/(\sqrt n u)\right)\nonumber\\[-1ex]
&&\\
&=& \Gamma(\alpha n) u^{\alpha n - 1} \sum^\infty_{m=1} b_{\alpha,n,m} \sin\left( \frac {m\pi}{\sqrt n u}\right) \nonumber .
\end{eqnarray}
Besides, an orthogonal expansion with Legendre polynomials for the density of $U^2$ is obtained by the moments
\begin{equation}
\label{eq:3b}   
\cE(U^{2k}) = \frac{\Gamma(\alpha n) k!}{\Gamma (\alpha n + 2k)} \sum_{(k)} \prod^{n}_{i=1} \frac{\Gamma (\alpha + 2k_i)}{\Gamma (\alpha) k_i!}\ ,
\end{equation}
derived from (\ref{eq:1b}).

Direct use of these orthogonal series within (\ref{eq:1e}) is possible to get the cdf of $S$, but the double series in theorem~\ref{theorem:4.1} is numerically more favourable. As mentioned in the introduction it would be theoretically more satisfying to have a single alternating series or a probability mixture for this cdf. A single power series is not directly available by a power series expansion of $G_{\alpha n}$ in (\ref{eq:1e}) since the moments of $\cot \Phi$ do not all exist. However, at least for integer $\alpha$, this obstacle is avoided by splitting the domain $[0, \arctan ( \sqrt{n-1})]$ of $\Phi$ by $\varphi_n = \arctan (1 / \sqrt{n-1})$, which is the maximal angle $\varphi$ for which the whole cone $\{\Phi \leq \varphi\}$ is contained completely within $\{x_1,\ldots,x_n | x_1,\ldots, x_n \geq 0\}$. For integer $\alpha$ the theorem~\ref{theorem:3.1} below provides a polynomial representation of the cdf of $\tan \Phi$ restricted to $0 \leq \tan \varphi \leq (n-1)^{-1/2}$. This simple representation is used within (\ref{eq:1e}) to integrate over $\tan \varphi \leq (n-1)^{-1/2}$. The integration over $\tan \varphi > (n-1)^{-1/2}$ needs only truncated moments of $\cot \Phi = (1- \cos^2\Phi)^{-1/2} \cos \Phi$, obtained by the binomial expansion with the truncated moments of $\cos \Phi$, given below in (\ref{eq:3f}).

\begin{theorem}
\label{theorem:3.1}
Let $\alpha$ be a positive integer, then
\begin{eqnarray}
\label{eq:3c}   
W_{\alpha,n}(t) &:=& \Pr \{ \tan \Phi \leq t\}\nonumber\\
&&\\
&=& \frac {\Gamma(\alpha n)}{(\Gamma(\alpha))^n n^{\alpha n/2}} \sum^{[(\alpha-1)n/2]}_{j=0} a_{\alpha,n,2j} \frac{t^{n-1+2j}}{n-1+2j}\ , \ 0\leq t \leq (n-1)^{-1/2},\nonumber
\end{eqnarray}
where the coefficients $a_{\alpha,n,2j}$ are given as the unique solutions of the linear equations
\begin{eqnarray}
\label{eq:3d}    
\lefteqn{\sum^{[(\alpha-1)n/2]}_{j=0} \Gamma \left( \frac 12 + \beta n + m - j\right) \Gamma\left( \frac{n-1}2 + j\right) a_{\alpha,n,2j}=}\nonumber\\
&&\\
&& \left\{ \begin{array}{ll}
\displaystyle 2 \frac {(2m)!}{n^m} \sum_{(m)} \prod^{n}_{j=1} \Gamma \left( 1/2 + \beta + m_j\right) / (2m_j)!\ , & \alpha = 2\beta + 1\\[5ex]
\displaystyle 2 \frac {(2m+n)!}{n^{m+n/2}} \sum_{(m)} \prod^{n}_{j=1} \Gamma \left( 1/2 + \beta + m_j\right) / (2m_j+1)!\ , & \alpha = 2\beta 
\end{array} \right\}\nonumber \\
&& m,j = 0 , \ldots, [(\alpha-1)n/2]. \nonumber
\end{eqnarray}
\end{theorem}
With any $\rho > 0$ the sums $\sum_{(m)}$ within the right hand sides of (\ref{eq:3d}) are also given by
\begin{eqnarray}
\label{eq:3e}     
\frac 1{\rho^m\pi} \int^{\pi}_{0} {\mathrm{Re}} \left\{ \left( \Gamma \left( \frac \alpha 2\right) M\left( \frac \alpha 2, \frac 12; \frac 14 \rho e^{i\varphi}\right)\right)^n e^{-im\varphi}\right\}d\varphi, \quad \alpha = 2\beta + 1, \nonumber\\
&&\\
\frac 1{\rho^m\pi} \int^{\pi}_{0} {\mathrm{Re}} \left\{ \left( \Gamma \left( \frac {\alpha+1} 2\right) M\left( \frac {\alpha+1} 2, \frac 32; \frac 14 \rho e^{i\varphi}\right)\right)^n e^{-im\varphi}\right\}d\varphi, \ \alpha = 2\beta , \nonumber
\end{eqnarray}
After Kummer's transformation (A.S. 13.1.27) the above confluent hypergeometric functions $M$ can also be expressed by Hermite polynomials due to (A.S. 13.6.17/18).

The truncated moments $\bar{\gamma}_{\alpha,n,k }=  \bar{\cE}((\cos \Phi)^{\alpha n + k})$, defined by integration over \mbox{$\varphi \geq \varphi_n$,} follow from theorem \ref{theorem:3.1} by straightforward calculation.
\begin{equation}
\label{eq:3f}    
\bar{\gamma}_{\alpha,n,k} = \gamma_{\alpha,n,k} - \frac{\Gamma(\alpha n)}{(\Gamma(\alpha))^n n^{\alpha n/2}} \frac 12 
\sum^{[(\alpha -1)n/2]}_{j=0} a_{\alpha,n,2j} B \left({\textstyle \frac{n-1}{2} +j, \frac{(\alpha-1) n+k+1}{2}-j; \frac{1}{n}}\right)
\end{equation}
with the incomplete beta function $B(a,b; x) = \int^{x}_{0} t^{a-1} (1-t)^{b-1} dt$. In particular with $\alpha = 1$ we have 
\begin{equation}
\label{eq:3g}    
\bar{\gamma}_{1,n,k} = \gamma_{1,n,k} - \frac{(n-1)!}{n^{n/2}} \frac{\pi^{(n-1)/2}}{\Gamma(\frac{n-1}2)} B \left( \frac{n-1}{2} , \frac{k+1}{2}; \frac{1}{n}\right).
\end{equation}

\paragraph{Two further remarks:} 
The equations (\ref{eq:3d}) can also be solved by explicit matrix inversion. The matrix $\big( \Gamma \left(\frac 12 + \beta n + m-j\right)\big)$, $m,j = 0, \ldots, N = [(\alpha -1) n/2]$, has the structure $\Gamma(x) {\mathrm{Diag}}\big(\ldots,(x)_m,\ldots\big)\big((x+m)_{N-j}\big)$
with factorials $(x)_m = x(x+1) \ldots (x+m-1)$. The inverse $(p_{mj}(x))$ of $((x+m)_{N-j})$ contains only polynomials $p_{mj}$ of degree $m$ in its $m^{\mathrm{th}}$ row. These polynomials can be obtained by interpolation from the $m+1$ values
\begin{eqnarray}
\begin{array}[t]{l} p_{mj}(-(N-k)) =\\[1ex]
(k = 0,\ldots, m) 
\end{array} 
\left\{
\begin{array}{clc}
\displaystyle \frac{(-1)^{j+k+1}}{(j+k-m)!(N-(j+k))!}& , & m-k\leq j \leq N-k\\[3ex]
0 & , & {\mathrm{otherwise}}
\end{array}
\right\}.\nonumber
\end{eqnarray}
The cdf $F_{1,n}$ of $U_{1,n}$ can be applied to a goodness of fit test of the hypothesis \mbox{$H_0: F = F_0$} against $H_1 : F \neq F_0$, where $F_0$ is any specified continuous cdf. Under $H_0$ the test statistic $\sum^{n}_{i=1}(F_0(X_{i:k}) - F_0(X_{i-1:k}))^2$, obtained from an ordered sample $X_{i:k}$, $i = 1, \ldots, k = n-1$, $(X_{0:k} = - \infty, \ X_{n:k} = \infty)$, has the same distribution as $U^2_{1,n}$. The asymptotic normality of $U_{1,n}^2$ can be derived from the asymptotic bivariate normal distribution of \mbox{$\big((Y_n-n)/\sqrt n$,} \mbox{$(Z_n-2n)/\sqrt{2n}\big)$} with correlation $\rho = 2 /\sqrt 5$. The convergence to a normal distribution is slow and might be accelerated  by a suitable transformation.

\section{The distribution of the sample variance}
\label{sec:4}
With the functions 
\begin{eqnarray}
\label{eq:4a}   
H_{\alpha,n,m}(r) &:=& \cE(G_{\alpha n +m} (r \sqrt n \cos \Phi))\nonumber \\
&&\\
&=& \frac {r^{\alpha n + m}}{\Gamma(\alpha n + m)} \sum^{\infty}_{k=0} \frac{\mu_{\alpha , n, m+k}}{\alpha n + m + k}\frac {(-r)^k}{k!}\ , \nonumber
\end{eqnarray}
$m \in \bN_0$, $\mu_{\alpha,n,k} = \cE(\sqrt n \cos \phi)^{\alpha n + k} = n^{(\alpha n + k)/2} \gamma_{\alpha,n,k}$ from (\ref{eq:2f}) we obtain
\begin{theorem}
\label{theorem:4.1}
Let $X_1,\ldots,X_n$ be independent random variables with density $(\Gamma(\alpha))^{-1} e^{-x} x^{\alpha -1}$, $\alpha > 0$, then the cdf of the sample variance $S^2 = \linebreak (n-1)^{-1} \sum^n_{i=1} (X_i - \bar{X})^2$ is given by
\begin{eqnarray}
\label{eq:4b}  
\lefteqn{\Pr\left\{ (n-1) S^2 \leq z = r^2\right\} = \sum^\infty_{j=0} \frac{\Gamma(\alpha n + 2j)}{\Gamma(\alpha n ) j! n^j} \left(\frac d{dz}\right)^j H_{\alpha,n,2j} (\sqrt z)}\\
\label{eq:4c}   
&=& \frac{r^{\alpha n}}{\Gamma(\alpha n)}\sum^\infty_{j=0} \left(\sum^\infty_{k=0} 
\begin{pmatrix}(\alpha n + k)/2 + j-1 \\ j \end{pmatrix}  
\frac{\mu_{\alpha,n,2j+k}}{\alpha n + k} \frac{(-r)^k}{k!}\right) \frac 1{n^j}\\
\label{eq:4d}    
&=& \sum^\infty_{j=0} \left(\sum^\infty_{k=0} \delta_{\alpha,n,j,k,\lambda} 
\begin{pmatrix}\alpha n + k -1 \\ k \end{pmatrix}
G_{\alpha n +k}  (\lambda r)\right) \left(\frac{\lambda^2}n\right)^j
\end{eqnarray}
with the alternating $k^{\mathrm{th}}$ order differences 
\begin{equation}
\label{eq:4e}  
\delta_{\alpha,n,j,k,\lambda} = \sum^k_{\ell =1} (-1)^\ell 
\begin{pmatrix}k \\ \ell\end{pmatrix}
\begin{pmatrix}(\alpha n + \ell)/2 + j -1 \\ j \end{pmatrix}
\frac{\mu_{\alpha,n,2j+\ell}}{\lambda^{\alpha n + 2j + \ell}}
\end{equation}
and any $\lambda > 0$. The series in (\ref{eq:4b}) is absolutely convergent with terms of the order $O \left(j^{-(n+1)/2}\right)$, where the $O$--constant depends only on $\alpha $ and $n$.  
\end{theorem}

As a simple numerical example we obtain with $\alpha = 1$ and $n = 10$ the value \mbox{$\Pr \{S \leq 2\} = 0.98530379\ldots$} with at least eight correct digits.

Finally, the method, explained before theorem~\ref{theorem:3.1}, is applied to obtain the alternative representations in theorem~\ref{theorem:4.2} for integer $\alpha$, including in particular a probability mixture for the exponential case $\alpha =1$. The degree of the polynomials in theorem~\ref{theorem:3.1} and the rate of convergence of the binomial series in (\ref{eq:4f}) limit the numerical use of theorem~\ref{theorem:4.2} to moderate values of $\alpha n$. 

The truncated moments
\begin{equation}
\label{eq:4f}  
\bar{M}_{\alpha,n,k} = \bar{\cE}\left( (\cot \Phi)^{\alpha n + k}\right) = \sum^\infty_{j=0} 
\begin{pmatrix} (\alpha n + k)/2+j-1 \\ j\end{pmatrix}
\bar{\gamma}_{\alpha,n,k+2j}
\end{equation}
are used with $\bar{\gamma}_{\alpha,n,k}$ from (\ref{eq:3f}), obtained by integration over $\varphi \geq \varphi_n = {\mathrm{arccot}} (\sqrt{n-1})$.

\begin{theorem}
\label{theorem:4.2}
With the truncated moments from (\ref{eq:4f}) and (\ref{eq:3f}), the pdf $w_{\alpha,n}$ and the cdf $W_{\alpha,n}$ of $\tan \Phi$ on $0 \leq \varphi \leq \varphi_n$ from theorem \ref{theorem:3.1}, the distribution of $S$ is given by each of the following three formulas (\ref{eq:4g}), (\ref{eq:4h}), (\ref{eq:4i}):
\begin{eqnarray}
\label{eq:4g}    
\lefteqn{\Pr\left\{ \sqrt{n-1} S \leq r\right\}}\nonumber\\
&&\\
&=& \frac{(r\sqrt n)^{\alpha n}}{\Gamma(\alpha n)} \sum^{\infty}_{k=0}(-1)^k \frac{\bar{\gamma}_{\alpha,n,k}}{\alpha n +k} P_{\alpha n,k} (r \sqrt n) + \int^{1/\sqrt{n-1}}_{0} G_{\alpha n} (r \sqrt n/t) w_{\alpha,n}(t)dt\nonumber
\end{eqnarray}
with the polynomials
\begin{eqnarray}
\label{eq:4h}    
\lefteqn{P_{\alpha n,k}(x) = \sum^{[k/2]}_{j=0} \begin{pmatrix} (\alpha n + k)/2 \\ j \end{pmatrix} 
\frac{x^{k-2j}}{(k-2j)!},} \nonumber\\[2ex]
\lefteqn{ \frac{(r \sqrt n)^{\alpha n}}{\Gamma (\alpha n)} \sum^{\infty}_{k=0} \frac{\bar{M}_{\alpha,n,k}}{\alpha n + k} \frac{(-r\sqrt n)^k}{k!}}\\
&+& W_{\alpha,n} (1/\sqrt{n-1}) G_{\alpha n} \left(r \sqrt{n(n-1)}\right) + \int^{\infty}_{r\sqrt{n(n-1)}} g_{\alpha n}(x) W_{\alpha,n} \left(r \sqrt n/x\right)dx . \nonumber
\end{eqnarray}
The infinite series in (\ref{eq:4h}) is also given by
\begin{equation}
\label{eq:4i}   
\sum^{\infty}_{k=0} \bar{\Delta}_{\alpha,n,k,\lambda} \begin{pmatrix} \alpha n + k -1 \\ k \end{pmatrix} G_{\alpha n + k}(\lambda r)
\end{equation}
with any $\lambda > 0$ and the alternating $k^{\mathrm{th}}$ order differences 
\begin{eqnarray}
 \bar{\Delta}_{\alpha,n,k,\lambda} &=& (-\Delta)^k \bar{M}_{\alpha,n,0,\lambda} = \sum^{k}_{\ell=0} (-1)^\ell \begin{pmatrix}k \\ \ell \end{pmatrix} \bar{M}_{\alpha,n,\ell,\lambda} \ , \nonumber\\
 \bar{M}_{\alpha,n,\ell,\lambda} &=& \bar{\cE} \left(\lambda^{-1} \sqrt n \cot \Phi\right)^{\alpha n + \ell}\ . \nonumber
\end{eqnarray}\\
\end{theorem}

\paragraph{Remarks:}
In particular, with $\alpha = 1$ the only coefficient in $w_{1,n}(t)$ is $a_{1,n,0} = \linebreak (n-1)b_{n-1}$ with the volume $b_{n-1} = \pi^{(n-1)/2}/\Gamma\left(\frac{n+1}2\right)$ of the $(n-1)$--unit ball. Then the last two terms in (\ref{eq:4h}) are reduced to
\begin{equation}
\label{eq:4j}    
\frac{(n-1)!}{\sqrt n(n(n-1))^{(n-1)/2}} b_{n-1} G_{n-1} (r\sqrt{n(n-1)}).
\end{equation}
The truncated moments $\bar{\gamma}_{1,n,k}$ were given in (\ref{eq:3g}). Only positive coefficients $\bar{\Delta}_{\alpha,n,k,\lambda}$ arise in (\ref{eq:4i}) with $\lambda = \sqrt{n(n-1)}$ since $1 > 1 - \frac{\cot \Phi}{\sqrt{n-1}} \geq 0$. 
The last term in (\ref{eq:4h}) is bounded by $W_{\alpha,n}(1/\sqrt{n-1})\left( 1-G_{\alpha n}(r\sqrt{n(n-1)})\right)$, which is often neglectible.  For actual computations the choice
\[
\lambda = \sqrt n\left( \frac{\bar\cE\left((\cot \Phi)^{\alpha n}\right)}{1 - W_{\alpha ,n}(1/\sqrt{n-1})}\right)^{1/(\alpha n)}
\]
might be more favourable. With this $\lambda$ the leading term in (\ref{eq:4i}) becomes $\left(1-W_{\alpha,n} \left( 1/\sqrt{n-1}\right)\right) G_{\alpha n} (\lambda r)$.

It would be desirable to find a similar probability mixture representation for the cdf of $S^2$ for all $\alpha > 0$ but presumedly, also in case of its existence, the computation of the required coefficients would not be easy.

\paragraph{Acknowledgement: }
 The author would like to thank S. Kotz for his  strong interest in this work and some hints at the literature within a personal communication.

\newcounter{appeqn}
\setcounter{appeqn}{1}
\newcounter{applem}
\setcounter{applem}{1}

\begin{appendix}
\section*{Appendix}

\renewcommand{\theequation}{\Alph{appeqn}.\arabic{equation}}
\setcounter{equation}{0}
\renewcommand{\thelemma}{\Alph{applem}.\arabic{lemma}}
\setcounter{lemma}{0}

\paragraph{\it {Proof of  theorem~\ref{theorem:3.1}}.} \ 
Let $\alpha$ be any positive integer and $\omega$ the Lebesgue measure on the $(n-1)$--unit sphere $S_{n-1}$. Then
\begin{eqnarray}
\label{eq:A.1}
\lefteqn{\int_{S_{n-1}} (x_1+\ldots+x_k)^k \left(\prod^{n}_{j=1} x_j\right)^{\alpha-1} d\omega}\nonumber\\
&&\\
&=& 2\Gamma((\alpha n + k)/2)^{-1} k! \sum_{(k)} \prod^{n}_{j=1} \Gamma (\alpha + k_j)/2)/k_j! \ , \quad \alpha + k_j \  \mathrm{odd}.\nonumber
\end{eqnarray}
After an orthogonal transformation $x = Ty$ with $y_1 = n^{-1/2} \sum^{n}_{i=1} x_i = \cos \vartheta_1$, followed by transformation to polar coordinates, we obtain with $c_j = \cos \vartheta_j$, $s_j = \sin \vartheta_j$ the same integral value as 
\begin{eqnarray}
\label{eq:A.2}
\lefteqn{\int_{S_{n-1}}(\sqrt n c_1)^k \prod^{n}_{j=1} (c_1/\sqrt n + s_1 \tau_j(\vartheta_2,\ldots,\vartheta_{n-1}))^{\alpha-1} d\omega}\nonumber\\
&=& n^{k/2} \int_{S_{n-1}} \sum^{(\alpha-1)n}_{\ell=0} \left( \sum_{(\ell)}\prod^{n}_{j=1}\begin{pmatrix} \alpha -1 \\ \ell_j\end{pmatrix} \tau^{\ell_j}_j\right) n^{(\ell-(\alpha-1)n)/2} \nonumber \\
&& \qquad \cdot \left( \prod^{n-1}_{j=2} s^{n-1-j}_j d\vartheta_j\right) c^{(\alpha-1)n+k-\ell}_1 s^{n-2+\ell}_1 d\vartheta_1 \nonumber\\
&=& n^{k/2} \sum^{[(\alpha-1)n/2]}_{j=0} a_{\alpha,n,2j} B\left( \frac{n-1}{2} + j, \frac{(\alpha-1) n+k+1}{2} -j\right). 
\end{eqnarray}
By comparison with (\ref{eq:A.1}) it follows
\begin{eqnarray}
\label{eq:A.3}
\lefteqn{\sum^{[(\alpha-1)n/2]}_{j=0} \Gamma \left( \frac{(\alpha-1)n + k +1}{2} -j\right) \Gamma \left(\frac{n-1}{2} + j\right) a_{\alpha,n,2j} } \nonumber\\
&&\\
&& \mbox{\hspace{2cm}} = 2n^{-k/2} k! \sum_{(k)} \prod^{n}_{j=1} \Gamma ((\alpha + k_j)/2)/k_j! \ . \nonumber
\end{eqnarray}
With $ \alpha = 2 \beta +1$ only even $k_j = 2m_j$ occur with sum $k = 2m$, and with $\alpha = 2\beta$ only $k_j = 2m_j + 1$ with sum $k = 2m+n$, which leads to the linear equations in (\ref{eq:3d}).

If $\ell$ is odd then the coefficients $a_{\alpha,n,\ell}$ vanish since
\[
\int_{S_{n-2}} \left( \prod^{n}_{j=1} \tau^{\ell_j}_j\right) \prod^{n-1}_{j=2} s_j^{n-1-j} d\vartheta_j
\]
is a linear combination of the integrals
\[
\int_{S_{n-2}} \left( \prod^{n-1}_{j=2}c^{m_j}_{j}\right)  s^{m_n}_{n-1} \prod^{n-1}_{j=2} s^{M_j +n - 1 - j}_{j} d\vartheta_j\ , \quad (M_{n-1}=0)\ ,
\]
which are different from zero only for even $m_2,\ldots,m_n$, $\sum^{n}_{j=2} m_j = \sum^{n}_{j=1} \ell_j = \ell$.

If $t = \tan \varphi \leq (n-1)^{-1}$ we obtain with $S_{n-1,\varphi} = S_{n-1} \cap \{ \vartheta_1 \leq \varphi\}$ that 
\begin{eqnarray*}
\Pr\{\tan \Phi \leq t\} &=& (\Gamma(\alpha))^{-n}  \int_{ \{\vartheta_1 \leq \varphi\}} \prod^{n}_{j=1} e^{-x_j} x^{\alpha-1}_j dx_j\\
&=& (\Gamma(\alpha))^{-n} \int^\infty_0 \left( \int_{S_{n-1},\varphi}\rho^{\alpha n-1} \exp(-\rho \sqrt n c_1) \prod^{n}_{j=1} (c_1/\sqrt n + s_1 \tau_j)^{\alpha-1} d\omega\right) d\rho \\
&=& \frac{\Gamma (\alpha n)}{(\Gamma(\alpha))^n n^{\alpha n /2}}  \int_{S_{n-1},\varphi} c^{-\alpha n}_1 \prod^{n}_{j=1}(c_1/\sqrt n + s_1 \tau_j)^{\alpha-1} d\omega\\
&=& \frac{\Gamma (\alpha n)}{(\Gamma(\alpha))^n n^{\alpha n /2}} \sum^{[(\alpha-1)n/2]}_{j=0} a_{\alpha,n,2j} \int^\varphi_0 s_1^{n-2+2j} c_1^{-n-2j} d\vartheta_1 \ ,
\end{eqnarray*}
where the integrals are given by $ (\tan \varphi)^{n-1+2j}/(n-1+2j$, which provides the asserted result in (\ref{eq:3c}). \hfill $\Box$\\

Comparing the series
\[
\sum^{\infty}_{k=0} \frac{(-1)^k}{k!}\frac{(z/y)^{(\alpha n+k)/2}}{\alpha n +k} = \Gamma (\alpha n) G_{\alpha n} (\sqrt{z/y})
\]
with the corresponding one for $\Gamma(\alpha n + 2j) G_{\alpha n +2j}(\sqrt{z/y})$ we obtain the identity
\begin{equation}
\label{eq:A.4}
\left(- \frac{\partial}{\partial y}\right)^j G_{\alpha n} \left( \sqrt{z/y}\right) = \frac{\Gamma(\alpha n + 2j)}{\Gamma(\alpha n)} \left( \frac{\partial}{\partial z}\right)^j G_{\alpha n +2j} \left( \sqrt{z/y}\right)
\end{equation}
by differentiation.

\begin{lemma}
\label{lemma:1}
\begin{equation}
\label{eq:A.5}
\left| \frac 1{j!} \left(- \frac{\partial}{\partial y}\right)^j G_{\alpha n} \left( \sqrt{z/y}\right)\right| \leq \frac{2^{\alpha n/2}}{\pi} \frac 1{jy^j}\ , \quad y,z >0 .
\end{equation}
\end{lemma}

\paragraph{\it Proof.}
Using (\ref{eq:A.4}) and Cauchy's integral formula for the derivatives of 
\[
\frac d{dz} \Gamma(\alpha n + 2j) G_{\alpha n + 2j} (\sqrt z) = \frac 12 z^{\alpha n/2 + j-1} \exp (-\sqrt z)
\]
the bound in (\ref{eq:A.5}) follows from 

\begin{eqnarray*}
\lefteqn{\Big| \frac 1{(j-1)!} \left( \frac d {dz}\right)^j \Gamma(\alpha n + 2j) G_{\alpha n + 2j} (\sqrt z)\Big|}\\
&=& \lim_{R\to \infty} \Big| \frac 1{4 \pi i} \oint_{w_R} \zeta^{\alpha n/2 +j-1} \exp(-\sqrt \zeta)(\zeta-z)^{-j} d\zeta \Big| \\
&\leq & \frac 1{2\pi} \int_0^\infty y^{\alpha n /2-1}\exp (-\sqrt{y/2}) dy = 2^{\alpha n/2} \Gamma(\alpha n)/ \pi\ ,
\end{eqnarray*}
where the way of integration is
\[
w_R = \left\{ {\mathrm{Re}}^{i\varphi} \big| - \frac \pi 2 \leq \varphi \leq \frac \pi 2\right\} \cup \left\{ iy \big| R \geq y \geq -R\right\}. \hspace{2cm} \Box
\]

\begin{lemma}
\label{lemma:2}
Let $\Phi$ denote the angle between $(X_1,\ldots,X_n)$ and $(1,\ldots,1)$ from (\ref{eq:1d}), $W_{\alpha,n}(t)$ the cdf of $\tan \Phi$ and $b_{n-1} = \pi^{(n-1)/2}/\Gamma(\frac {n+1}2)$ the volume of the $(n-1)$--unit ball, then 
\begin{equation}
\label{eq:A.6}
W_{\alpha,n}(t) \sim \frac {\Gamma(\alpha n) b_{n-1}} {(\Gamma(\alpha))^n n^{(\alpha - 1/2)n}} t^{n-1}, \quad t \downarrow 0,
\end{equation}
and in particular 
\begin{equation}
\label{eq:A.7}
W_{1,n}(t) = \frac{(n-1)!}{n^{n/2}} b_{n-1} t^{n-1}, \quad 0 \leq t \leq (n-1)^{-1/2} .
\end{equation}
Furthermore
\begin{equation}
\label{eq:A.8}
\cE \left((\cos \Phi)^s\right) \sim \frac{\Gamma(\alpha n)(2\pi)^{(n-1)/2}}{(\Gamma(\alpha))^n n^{(\alpha -1/2)n}} s^{-(n-1)/2}, \quad s\uparrow \infty. 
\end{equation}
\end{lemma}

\paragraph{\it Proof.}
Let $x,x_g$ denote the arithmetic and the geometric mean of $(x_1,\ldots,x_n \geq 0$ and $\varphi$ the angle between $ (x_1,\ldots,x_n)$ and $(1,\ldots,1)$. Setting $x_i = (1+\epsilon_i)x$ with any fixed $x$ we have the relation 
\[
x^{-2} \sum^{n}_{i=1} (x_i -x)^2 = \sum^{n}_{i=1} \epsilon^2_i = n \tan^2 \varphi
\]
and for $\varphi \downarrow 0$
\[
x^n_g = x^n \prod^n_{i=1} (1+\epsilon_i) \sim x^n \exp\left( \frac 12 \sum^n_{i=1} \epsilon^2_i\right) \sim x^n.
\]
It follows
\begin{eqnarray*}
W_{\alpha,n}(t) &=& \Pr \{\Phi \leq \varphi\} = \int_{\Phi\leq \varphi} \prod^n_{i=1} g_\alpha(x_i) dx_i\\
&=& (\Gamma(\alpha))^{-n}\int_{\Phi\leq \varphi} x^{n(\alpha-1)}_g \exp (-nx) dx_1\ldots dx_n\\
&\sim& (\Gamma(\alpha))^{-n} b_{n-1} \sqrt n \int_0^\infty x^{(\alpha-1)n} \exp(-nx) (x\sqrt n \tan \varphi)^{n-1} dx
\end{eqnarray*}
which provides (\ref{eq:A.6}).

For $\alpha =1$ this asymptotic relation can be replaced by an equation if $\varphi \leq \varphi_n = \arctan (1 /\sqrt{n-1})$, which is the largest angle $\varphi$ for which the whole cone $\{\Phi \leq \varphi\}$ is contained within \mbox{$\{x_1,\ldots,x_n|x_1,\ldots,x_n \geq 0\}$.}

To prove (\ref{eq:A.8}) we obtain from (\ref{eq:A.6})
\begin{eqnarray*}
\Pr\{|\ln (\cos \Phi)| \leq \epsilon\} &=& \Pr \{\ln (1+\tan^2 \Phi) \leq 2\epsilon\}\\
&\sim & W_{\alpha n} \left(\sqrt{2\epsilon}\right) \sim c_{\alpha,n} \epsilon^{(n-1)/2}, \quad \epsilon \downarrow 0,
\end{eqnarray*}
with a factor $c_{\alpha,n}$ determined by (\ref{eq:A.6}). For $s \to \infty$ the relation (\ref{eq:A.8}) follows from the Hardy--Littlewood--Karamata Tauber theorem (cf. e.g. chapter 13 in Feller 1971), applied to
\[
\cE ((\cos \Phi)^s) = \cE\left( \exp(-s|\ln(\cos \Phi)|)\right). \hspace{4cm} \Box
\] 

\paragraph{\it Proof of theorem \ref{theorem:4.1}. }  
With any $\vartheta \in (0,1)$ \quad $\cE_U \left( G_{\alpha n}\big((z/(U^2-\vartheta/n))^{1/2}\big)\right)$ can be represented by the Taylor series
\[
\sum^\infty_{j=0} \frac 1{j!} \cE_U\left( \left(- \frac \partial {\partial y}\right)^j G_{\alpha n} \left( \sqrt{z/y}\right)\Big|_{y=U^2} \right) \left(\frac\vartheta n\right)^j.
\]
With $U^{-1} = \sqrt n \cos \Phi$ and the bounds from Lemma~\ref{lemma:1} the terms of this series are absolutely bounded by
\[
\frac{2^{\alpha n/2}}{\pi}\cE\left( (\cos \Phi)^{2j}\right) \vartheta^j/j = O (\vartheta^j j^{-(n+1)/2})
\]
due to (\ref{eq:A.8}) with $s = 2j$ and on $O$--constant depending only on $\alpha$ and $n$.

For $\vartheta \uparrow 1$ it follows from (\ref{eq:1f}) and (\ref{eq:A.4}) that 
\begin{eqnarray*}
\Pr\{(n-1)S^2 \leq z = r^2\} 
&=& \cE _U\left( G_{\alpha n} \left( \sqrt{\frac z{U^2 - 1/n}} \right)\right)\\
&=& \sum^\infty_{j=0} \frac{\Gamma(\alpha n + 2j)}{\Gamma(\alpha n) j! n^j}\cE_U \left(\left( \frac \partial {\partial z}\right)^j G_{\alpha n + 2j} \left( \sqrt z/U\right)\right)\\
&=& \sum^\infty_{j=0} \frac{\Gamma(\alpha n + 2j)}{\Gamma(\alpha n) j! n^j} \left( \frac d {dz}\right)^j H_{\alpha,n,2j} (\sqrt z)
\end{eqnarray*}
with the functions $H_{\alpha,n,2j}$ from (\ref{eq:4a}).

With any $\lambda > 0$ we obtain
\begin{eqnarray*}
\lefteqn{\frac d{dr} \left( \frac{\Gamma(\alpha n + 2j)}{\Gamma (\alpha n) j! n^j} \left( \frac d{dz} \right)^j H_{\alpha,n,2j} (\sqrt z) \Big|_{z=r^2}\right)}\\
&=& \frac \lambda {\Gamma(\alpha n)} \sum^\infty_{k=0} \begin{pmatrix}(\alpha n + k)/2+j-1 \\ j\end{pmatrix} 
\frac{\mu_{\alpha,n,2j + k}}{\lambda^{\alpha n + 2j +k}} \frac{(-1)^k}{k!}(\lambda r)^{\alpha n -1+k} \left(\lambda^2/n\right)^j.
\end{eqnarray*}
After multiplication by $\exp(-\lambda r) \sum^\infty_{\ell = 0} (\lambda r)^\ell/\ell !$ the series (\ref{eq:4d}) is obtained by integration over $r$. \hfill $\Box$

\paragraph{\it Proof of  theorem \ref{theorem:4.2}. } \ 
With the pdf $w_{\alpha,n}$ and the cdf $W_{\alpha,n}$ of $\tan \Phi$ from theorem~\ref{theorem:3.1} and the truncated moments $\bar{M}_{\alpha,n,k}= \bar{\cE} \left( (\cot \Phi)^{\alpha n + k}\right)$ from (\ref{eq:4f}) it follows with (\ref{eq:1e}) that
\begin{eqnarray*}
\Pr \{ \sqrt{n-1} S \leq r\} &=& \cE\left(G_{\alpha n} (r \sqrt n \cot \Phi)\right)\\
&=& \int^{1/\sqrt{n-1}}_0 G_{\alpha n} (r \sqrt n /t) w_{\alpha,n} (t) dt + \frac{(r\sqrt n)^{\alpha n}}{\Gamma(\alpha n)} \sum^\infty_{k=0} \frac{\bar{M}_{\alpha,n,k}}{\alpha n + k} \frac{(-r\sqrt n)^k}{k!}.
\end{eqnarray*}
After integration by parts and the substitution $t \to r\sqrt n /x$ the above integral becomes
\[
W_{\alpha n} \left( 1 /\sqrt{n-1}\right) G_{\alpha n} \left(r\sqrt{n(n-1)}\right) + \int^\infty_{r\sqrt{n(n-1)}} g_{\alpha n}(x) W_{\alpha n} \left( r\sqrt n /x\right) dx.
\]
Inserting the series for $\bar{M}_{\alpha,n,k}$ from (\ref{eq:4f}) formula (\ref{eq:4g}) follows by rearranging the resulting absolutely convergent double series according to the truncated moments $\bar{\gamma}_{\alpha,n,k} = \bar{\cE} \left( (\cos \Phi)^{\alpha n+k}\right)$ from (\ref{eq:3f}). 

The series in (\ref{eq:4i}) is obtained by the same way as (\ref{eq:2ma}). \hfill $\Box$
\end{appendix}

%
                    


\end{document}